\documentclass[graybox]{svmult}

\usepackage{mathrsfs}
\usepackage{bbm}
\usepackage{bm}

\usepackage{makeidx}         % allows index generation
\usepackage{graphicx}        % standard LaTeX graphics tool
\usepackage{multicol}        % used for the two-column index
\usepackage[bottom]{footmisc}% places footnotes at page bottom

\usepackage{newtxtext}       %
\usepackage{newtxmath}       % selects Times Roman as basic font

\makeindex             % used for the subject index
\usepackage{verbatim}
\usepackage[all,cmtip]{xy}
\usepackage[mathscr]{euscript}

\usepackage{etex}
\usepackage{graphicx}
\usepackage{pictex}
\usepackage{epstopdf}
\usepackage{color}
\date{\today}

\renewcommand{\D}{{\mathbb D}}
\newcommand{\N}{{\mathbb N}}
\newcommand{\Z}{{\mathbb Z}}

\renewcommand{\Re}{{\mathrm{Re} \,}}

\def\jdf#1{\textcolor{blue}{#1}} %comments by Jake Fillman
\def\szg#1{\textcolor{cyan}{#1}} %comments by Shuzheng

\spnewtheorem{thm}{Theorem}[section]{\bf }{}
\spnewtheorem{coro}[thm]{Corollary}{\bf }{}
\spnewtheorem{prop}[thm]{Proposition}{\bf }{}
\spnewtheorem{lem}[thm]{Lemma}{\bf }{}

\spnewtheorem*{HDC}{Hausdorff Dimension Conjecture}{\bf }{}

\allowdisplaybreaks

\newcommand{\supp}{{\mathrm{supp}}}

\newcommand{\set}[1]{\left\{#1\right\}}

\begin{document}

\title*{On Simon's Hausdorff Dimension Conjecture}

\author{David Damanik \and Jake Fillman \and Shuzheng Guo \and Darren C.\ Ong}
 \authorrunning{D.\ Damanik \and J.\ Fillman \and S.\ Guo \and D.C.\ Ong}

\institute{David Damanik
\at  Rice University, Houston, TX~77005, USA
\at \email{damanik@rice.edu}
\and Jake Fillman
\at Texas State University, San Marcos, TX~78666, USA
\at \email{fillman@txstate.edu}
\and Shuzheng Guo
\at Ocean University of China, Qingdao 266100, Shandong, China and Rice University, Houston, TX~77005, USA
\at \email{gszouc@gmail.com}
\and Darren C.\ Ong
\at Xiamen University Malaysia, 43900 Sepang, Selangor Darul Ehsan, Malaysia
\at \email{darrenong@xmu.edu.my}}

\maketitle

%\author[D.\ Damanik]{David Damanik}
%\address{Department of Mathematics, Rice University, Houston, TX~77005, USA}
%\email{damanik@rice.edu}
%\thanks{D.D.\ was supported in part by NSF grant DMS--1700131 and by an Alexander von Humboldt Foundation research award}

%\author[J.\ Fillman]{Jake Fillman}
%\address{Department of Mathematics, Texas State University, San Marcos, TX~78666}
%\email{fillman@txstate.edu}
%\thanks{J.F.\ was supported in part by Simons Foundation Collaboration Grant \#711663}

%\author[S.\ Guo]{Shuzheng Guo}
%\address{Ocean University of China, Qingdao 266100, Shandong, China and Rice University, Houston, TX~77005, USA}
%\email{gszouc@gmail.com}
%\thanks{S.G.\ was supported by CSC (No. 201906330008) and NSFC (No. 11571327)}

%\author[D.\ Ong]{Darren C.\ Ong}
%\address{Department of Mathematics, Xiamen University Malaysia, 43900 Sepang, Selangor Darul Ehsan, Malaysia}
%\email{darrenong@xmu.edu.my}
%\thanks{D.O.\ was supported in part by a grant from the Fundamental Research Grant Scheme from the Malaysian Ministry of Education (Grant No: FRGS/1/2018/STG06/XMU/02/1) and a Xiamen University Malaysia Research Fund (Grant Number: XMUMRF/2020-C5/IMAT/0011)}

\abstract{Barry Simon conjectured in 2005 that the Szeg\H{o} matrices, associated with Verblunsky coefficients $\{\alpha_n\}_{n\in\Z_+}$ obeying $\sum_{n = 0}^\infty n^\gamma |\alpha_n|^2 < \infty$ for some $\gamma \in (0,1)$, are bounded for values $z \in \partial \D$ outside a set of Hausdorff dimension no more than $1 - \gamma$. Three of the authors recently proved this conjecture by employing a Pr\"ufer variable approach that is analogous to work Christian Remling did on Schr\"odinger operators. This paper is a companion piece that presents a simple proof of a weak version of Simon's conjecture that is in the spirit of a proof of a different conjecture of Simon.}

\section{Introduction}

This paper is concerned with Barry Simon's Hausdorff dimension conjecture for orthogonal polynomials on the unit circle (OPUC).

There has been a large amount of activity studying OPUC in the past two decades, primarily due to the monographs \cite{S04, S05} by Simon, to which we refer the reader for general background material.
In these monographs, Simon makes a number of conjectures, which are listed for the reader's convenience in \cite[Appendix~D]{S05}.
Our main purpose here is to prove a result related to one of these conjectures, the Hausdorff dimension conjecture; see \cite[Conjecture~D.3.5, p.~982]{S05}.

Before stating this conjecture, let us describe the general setting of OPUC.
Suppose $\mu$ is a non-trivial (i.e., not finitely supported) probability measure on the unit circle $\partial \mathbb{D} = \{ z \in \mathbb{C} : |z| = 1 \}$.
By the non-triviality assumption, the functions $1, z, z^2, \cdots$ are linearly independent in the Hilbert space $\mathcal{H} = L^2(\partial\mathbb{D}, \mu)$, and hence one can form, via the Gram-Schmidt procedure, the \emph{monic orthogonal polynomials} $\Phi_n(z)$, whose Szeg\H{o} dual is defined by
$$\Phi_n^{*} = z^n\overline{\Phi_n({1}/{\overline{z}})}.$$
There are constants $\{\alpha_n\}_{n\in\Z_+}$ in $\mathbb{D}=\{z\in\mathbb{C}:|z|<1\}$, called the \emph{Verblunsky coefficients}, so that
\begin{equation}\label{eq01}
\Phi_{n+1}(z) = z \Phi_n(z) - \overline{\alpha}_n \Phi_n^*(z), \qquad \textrm{ for } n\in \Z_+,
\end{equation}
which is the so-called \emph{Szeg\H{o} recurrence} (here and throughout the paper, $\Z_+$ denotes the set of all non-negative integers and $\N$ will denote the set of positive integers). Conversely, every sequence $\{\alpha_n\}_{n\in\Z_+}$ in $\mathbb{D}$ arises as the sequence of Verblunsky coefficients for a suitable nontrivial probability measure on $\partial\mathbb D$.

If we consider instead the orthonormal polynomials
$$\varphi_n(z)=\frac{\Phi_n(z)}{\|\Phi_n(z)\|_{\mu}},$$
where $\|\cdot\|_{\mu}$ is the norm of $\mathcal{H}$, one can verify that \eqref{eq01} becomes
\begin{equation}\label{eq01b}
\rho_n \varphi_{n+1}(z) = z \varphi_n (z) - \overline{\alpha}_n \varphi_n^*(z), \textrm{ for } n\in \Z_+,
\end{equation}
where $\rho_n = (1 - |\alpha_n|^2)^{1/2}$.

The Szeg\H{o} recurrence can be written in a matrix form as follows:
\begin{equation} \label{eq:SzegoVarphiMatRec}
\begin{bmatrix}
\varphi_{n+1}(z)\\
\varphi_{n+1}^{*}(z)
\end{bmatrix}
=
\frac{1}{\rho_n}
\begin{bmatrix}
z & -\overline{\alpha}_n\\
-\alpha_n z & 1
\end{bmatrix}
\begin{bmatrix}
\varphi_n(z)\\
\varphi_n^{*}(z)
\end{bmatrix}, \textrm{ for } n\in\Z_+.
\end{equation}

Alternatively, one can consider a different initial condition and derive the \emph{orthogonal polynomials of the second kind}, by setting $\psi_0(z) = 1$ and then
\begin{equation} \label{eq:SzegoPsiMatRec}
\begin{bmatrix}
\psi_{n+1}(z)\\
- \psi^{*}_{n+1}(z)
\end{bmatrix}
=
\frac{1}{\rho_n}
\begin{bmatrix}
z & -\overline{\alpha}_n\\
-\alpha_n z & 1
\end{bmatrix}
\begin{bmatrix}
\psi_n(z)\\
- \psi_n^{*}(z)
\end{bmatrix}, \textrm{ for } n\in\Z_+.
\end{equation}
In particular, the sequence $\{\psi_n(z)\}_{n \in \Z_+}$ is precisely the same as the first-kind polynomials for the measure $\widetilde\mu$ with Verblunsky coefficients $\widetilde\alpha_n =-\alpha_n$; compare the discussion in \cite[Section~3.2, p.~222]{S04}, especially Equations~(3.2.2) and (3.2.3).
Define
\begin{equation} \label{eq:TnzDef}
T_n(z)=\frac{1}{2}
\begin{bmatrix}
\varphi_n(z) + \psi_n(z) & \varphi_n(z) - \psi_n(z)\\
\varphi^*_n(z) - \psi^*_n(z) & \varphi^*_n(z) + \psi^*_n(z)
\end{bmatrix}, \textrm{ for } n\in\Z_+.
\end{equation}
%
%For a Borel set $A\subset [0,1]$, and $0 \le d \le 1$, define
%\[
%S_{d,\delta}(A) = \inf\left\{\sum_{j=0}^{\infty} |b_j|^d | A\subset \bigcup_j b_j, |b_j|<\delta\right\}.
%\]
%The inf is taken over all covers by intervals of size at most $\delta$. Obviously, this inf increases as $\delta$ decreases. Recall that the  \emph{$\alpha$-dimensional Hausdorrf measure} is definded by
%\[
%h_d(A) = \sup_{\delta} S_{d,\delta}(A).
%\]
%For any $A$, there is a unique number, $d_H(A)$, called the \emph{Hausdorff dimension} of $A$, so that
%\[
%h_d(A) = \begin{cases}
%0,\qquad & d>d_{H}(A),\\
%\infty, & d<d_{H}(A).
%\end{cases}
%\]
%$h_{d(A)}(A)$ can be $0$, finite, or $\infty$.

We can now state \cite[Conjecture~D.3.5, p.~982]{S05}:

\begin{HDC}
Assuming that the Verblunsky coefficients obey a decay estimate of the form
\begin{equation}\label{e.vcass}
\sum_{n = 0}^\infty n^\gamma |\alpha_n|^2 < \infty
\end{equation}
for some $\gamma \in (0,1)$, the associated Szeg\H{o} matrices $T_n(z)$ are bounded in $n \in \Z_+$ for values $z \in \partial \D$ outside a set of Hausdorff dimension no more than $1 - \gamma$.
\end{HDC}

The primary interest in such a statement comes from the general principle that says that the singular part of $\mu$ is supported by the set of those $z \in \partial \D$ for which the associated Szeg\H{o} matrices are unbounded, that is, such a result will imply that the singular part of the measure, $\mu_\mathrm{sing}$, has a support of Hausdorff dimension at most $1 - \gamma$. As discussed in \cite{S05}, the upper bound on the Hausdorff dimension of the exceptional set of spectral parameters for which the Szeg\H{o} matrices are unbounded is tight, that is, it cannot be improved.

Three of the authors recently proved Simon's Hausdorff dimension conjecture \cite{GDO20}. While the proof of the full conjecture (which is inspired by work of Remling \cite{r, r2}) is not short and somewhat technical, we show here that if one is willing to give up a factor $2$, then there is a rather short proof. %In fact, the proof can accommodate logarithmic divergence of the series in \eqref{e.vcass}.
In other words, we will show how to quickly establish the following result:

\begin{thm}\label{t.main}
Suppose that $\mu$ is such that the associated Verblunsky coefficients satisfy \eqref{e.vcass}. Then there is a set $S \subset \partial \D$ of Hausdorff dimension at most $2(1 - \gamma)$ so that for $z \in \partial \D \setminus S$,
$$
\sup_{n \ge 0} \|T_n(z)\| < \infty.
$$
In particular, $\mu_\mathrm{sing}$ is supported by a set of dimension at most $2(1 - \gamma)$.
\end{thm}

Of course, the theorem is only meaningful for $\gamma \in (1/2,1)$, so by giving up a factor of $2$ in the estimate, one not only loses the tightness of the estimate, one also reduces the size of the relevant range of $\gamma$'s by half. Nevertheless, we feel it is still worthwhile to point out that this weaker result can be obtained via a short argument, which avoids the lengthy and intricate arguments used in \cite{GDO20}. This short argument is inspired by the work \cite{d}, which incidentally proves another conjecture of Simon, namely \cite[Conjecture~D.3.4, p.~982]{S05}.
Earlier related work can be found in \cite{dk}.
We will prove Theorem~\ref{t.main} in Section~\ref{sec.2}.

Let us briefly remark that one can extend Theorem~\ref{t.main} to cover the case of logarithmic divergence.
\begin{coro} \label{coro:main}
Suppose that $\mu$ is such that the associated Verblunsky coefficients satisfy \begin{equation} \label{eq:vcasslogdiv} \sum_{n=0}^N n^\gamma |\alpha_n|^2 \leq A (\log N)^B \end{equation} for all $N \geq 2$, where $A,B>0$, $\gamma \in (0,1)$ are constants. Then there is a set $S \subset \partial \D$ of Hausdorff dimension at most $2(1 - \gamma)$ so that for $z \in \partial \D \setminus S$,
$$
\sup_{n \ge 0} \|T_n(z)\| < \infty.
$$
In particular, $\mu_\mathrm{sing}$ is supported by a set of dimension at most $2(1 - \gamma)$.
\end{coro}

All four of the authors of this manuscript have ties to Texas, so we are familiar with and grateful for Lance's contributions to advancing mathematics in this state in his role at Baylor University. Many happy returns, Lance!

\subsection*{Acknowledgements}D.D.\ was supported in part by NSF grant DMS--1700131 and by an Alexander von Humboldt Foundation research award. J.F.\ was supported in part by Simons Foundation Collaboration Grant \#711663. S.G.\ was supported by CSC (No. 201906330008) and NSFC (No. 11571327). D.O.\ was supported in part by two grants from the Fundamental Research Grant Scheme from the Malaysian Ministry of Education (Grant Numbers: FRGS/1/2018/STG06/XMU/02/1 and FRGS/1/2020/STG06/XMU/02/1) and a Xiamen University Malaysia Research Fund (Grant Number: XMUMRF/2020-C5/IMAT/0011).

\section{A Weak Version of Simon's Hausdorff Dimension Conjecture}\label{sec.2}

%%%%%%%%%%%%%%%%%%%%%%%%%%
\subsection{A Basic Estimate}

We start by deducing a basic consequence of the assumption \eqref{e.vcass}. Define $d:=1-\gamma$.

\begin{lem}\label{l.dyadicCSest}
Under the assumption \eqref{e.vcass}, $\{n^{-(d/2+\varepsilon/4)} \alpha_n\}_{n \in \N} \in \ell^1(\N)$ for all $\varepsilon > 0$.
\end{lem}

\begin{proof}
Applying the Cauchy-Schwarz inequality to dyadic blocks, for example, we see that
\begin{align*}
\sum_{n = 1}^\infty n^{-(d/2+\varepsilon/4)} |\alpha_n| & = \sum_{k = 0}^\infty \sum_{n = 2^k}^{2^{k+1}-1} n^{-(d/2+\varepsilon/4)} |\alpha_n| \\
& = \sum_{k = 0}^\infty \left( \sum_{n = 2^k}^{2^{k+1}-1} n^{-(1/2 + \varepsilon/4)} n^{\gamma/2} |\alpha_n| \right) \\
& \le \sum_{k = 0}^\infty \left( \sum_{n = 2^k}^{2^{k+1}-1} n^{-(1+\varepsilon/2)} \right)^{1/2} \left( \sum_{n = 2^k}^{2^{k+1}-1} n^{\gamma} |\alpha_n|^2 \right)^{1/2} \\
& \lesssim \sum_{k = 0}^\infty 2^{-k\varepsilon/4}  \\
& < \infty,
\end{align*}
as desired.
\end{proof}

%%%%%%%%%%%%%%%%%%%

\subsection{Pr\"ufer Variables}
%%%%%%%%%%%%%%%%%%%%%%%%%%%%%%%%%%%%%%%%%%%%%%%%%%%%%%%%%%%%%%%%%%%%%%%%%%%%%%%%%%%%%

Let $\{\alpha_n\}_{n\in\Z_+}$ be the Verblunsky coefficients of a nontrivial probability measure
$\mu$ on $\partial \D$. As mentioned above, the $\alpha$'s give rise to a sequence
$\{\Phi_n(z)\}_{n\in\Z_+}$ of monic polynomials (via the Szeg\H{o} recurrence) that are orthogonal
with respect to $\mu$.
For $\beta \in [0,2\pi)$, we also consider the monic polynomials
$\{\Phi_n(z,\beta)\}_{n\in\Z_+}$ that are associated in the same way with the rotated Verblunsky
coefficients $\{ e^{i\beta} \alpha_n \}_{n\in\Z_+}$.

Let $\eta \in [0,2\pi)$. Define the Pr\"ufer variables by
$$
\Phi_n(e^{i\eta},\beta) = R_n(\eta,\beta) \exp \left[ i(n \eta + \theta_n(\eta,\beta))
\right],
$$
where $R_n > 0$, $\theta_0 \in [0,2\pi)$, and $|\theta_{n+1} - \theta_n| < \pi/2$ (compare \cite[Corollary~10.12.2]{S05}).
These variables obey the following pair of equations:
\begin{align*}
\frac{R_{n+1}^2(\eta,\beta)}{R_n^2(\eta,\beta)} & = 1 + | \alpha_n |^2 - 2 \Re
\left(\alpha_n e^{i[(n+1)\eta + \beta + 2
\theta_n(\eta,\beta)]} \right), \\
e^{-i(\theta_{n+1}(\eta,\beta) - \theta_n(\eta,\beta))} & = \frac{1 - \alpha_n
e^{i[(n+1)\eta + \beta + 2 \theta_n(\eta,\beta)]}}{\left[1 + | \alpha_n |^2 - 2 \Re
\left(\alpha_n e^{i[(n+1)\eta + \beta + 2 \theta_n(\eta,\beta)]} \right)\right]^{1/2}}.
\end{align*}
We also define $r_n(\eta,\beta) = |\varphi_n(\eta,\beta)|$.

When $\{ \alpha_n \} \in \ell^2$,
\begin{equation}\label{fs}
r_n(\eta,\beta) \sim R_n(\eta,\beta) \sim \exp \left( - \sum_{j=0}^{n-1} \Re (\alpha_j
e^{i[(j+1)\eta + \beta + 2 \theta_j(\eta,\beta)]}) \right).
\end{equation}
(We write $f_n \sim g_n$ if there is $C > 1$ such that $C^{-1}g_n \le f_n \le C g_n$ for
all $n$.) For the Pr\"ufer equations and \eqref{fs}, see \cite[Theorems~10.12.1 and
10.12.3]{S05}.

%%%%%%%%%%%%%%%%%%%%%%%%%%%%%%%%%%%%%%%%%%%%%%%%%%%%%%%%%%%%%%%%%%%%%%%%%%%%%%%%%%%%%
\subsection{Unboundedness and Infinite Energy}
%%%%%%%%%%%%%%%%%%%%%%%%%%%%%%%%%%%%%%%%%%%%%%%%%%%%%%%%%%%%%%%%%%%%%%%%%%%%%%%%%%%%%

In this section, we will prove that the set of $\eta\in (0,2\pi)$ for which the radius is unbounded has Hausdorff dimension no more than $2d$, which is stated as follows. The overall strategy in our proof of this statement will be inspired by Damanik-Killip \cite{dk}.

\begin{prop}\label{p31}
Assume \eqref{e.vcass}. Then the set
$$
S = \{ \eta \in [0,2\pi) : R_n(\eta,\beta) \text{ is unbounded for some } \beta \}
$$
has Hausdorff dimension no more than $2d = 2(1-\gamma)$.
\end{prop}

By \eqref{e.vcass}, $\{\alpha_n\}_{n\in\Z_+} \in \ell^2$. Therefore, because of
\eqref{fs}, it suffices for our purposes to show that
$$
A(n,\eta,\beta) = \sum_{j=0}^{n-1} \alpha_j e^{i[(j+1)\eta + \beta + 2
\theta_j(\eta,\beta)]}
$$
is a bounded function of $n$ for all $\beta$, provided that $\eta$ is away from a set of
Hausdorff dimension at most $2d$.

\begin{lem}\label{klstool}
If
$$
\widehat{\alpha}(\eta,n) = \lim_{N \to \infty} \sum_{j = n}^N \alpha_j e^{i j \eta}
$$
exists and obeys
\begin{equation}\label{fsr}
\sum_{j=1}^\infty | \widehat{\alpha} (\eta,j) \alpha_{j-1} | < \infty,
\end{equation}
then $\eta \not\in S$.
\end{lem}

\begin{proof}
We will show that $A(n,\eta,\beta)$ is bounded (in $n$) for every $\beta \in [0,2\pi)$
when \eqref{fsr} holds. The assertion then follows from \eqref{fs}.

Write $\gamma_j(\eta,\beta) = (j+1)\eta + \beta + 2 \theta_j(\eta,\beta)$. We have
\begin{align*}
A(n,\eta,\beta) & = \sum_{j=0}^{n-1} \left[\widehat{\alpha} (\eta,j) -
\widehat{\alpha} (\eta,j+1)\right] e^{i\gamma_j(\eta,\beta) - i j \eta} \\
& = \sum_{j=1}^{n-1} \widehat{\alpha} (\eta,j) \left[ e^{i\gamma_j(\eta,\beta)} -
e^{i(\gamma_{j-1}(\eta,\beta) + \eta)} \right] e^{-i j \eta} + O(1).
\end{align*}
Since
\begin{align*}
| e^{i\gamma_j(\eta,\beta)} - e^{i(\gamma_{j-1}(\eta,\beta) + \eta)}| & \le |
\gamma_j(\eta,\beta) - \gamma_{j-1}(\eta,\beta) - \eta | \\ & = 2 | \theta_j(\eta,\beta)
- \theta_{j-1}(\eta,\beta) | \\ & \lesssim |\alpha_{j-1}|,
\end{align*}
where the first inequality follows from the mean value theorem and the last inequality follows from \cite[Corollary 10.12.2]{S05} as well as the fact that $\alpha$'s are uniformly bounded away from $1$, boundedness of $A(n,\eta,\beta)$ follows.
\end{proof}

\begin{lem}\label{SZ}
Let $\nu$ be a positive measure on $[0,2\pi)$. For each $s \in (0,1)$ and every
measurable function $m : [0,2\pi) \to \Z_+$,
$$
\Biggl\{ \int\ \Biggl| \sum_{n=0}^{m(\eta)}  c_n e^{-in\eta} \Biggr| \, d\nu(\eta)
\Biggr\}^2 \lesssim \mathcal{E}_s (\nu) \sum_{n=0}^\infty (n+1)^{1-s}
\big|c_n\big|^2 ,
$$
where $\mathcal{E}_s(\nu)$ denotes the $s$-energy of $\nu$, which is defined by
$$
\mathcal{E}_s (\nu) = \iint (1+|x-y|^{-s}) \, d\nu(x) \, d\nu(y).
$$
\end{lem}

\begin{proof}
This follows by slightly adjusting the calculation from \cite[\S XIII.11,
p.~196]{Zygmund} (see also \cite[\S V.5]{Carleson}).
\end{proof}

\begin{proof}[Proof of Proposition~\ref{p31}.]
We will apply the criterion of Lemma~\ref{klstool}. Let us first note that by the theorem
of Salem-Zygmund \cite[Theorem XIII.11.3(2)]{Zygmund} and the connection between capacity and Hausdorff measure \cite[\S IV.1]{Carleson}, the series defining $\widehat{\alpha}$ converges off a set of Hausdorff dimension $d$. Therefore, we may exclude from consideration those values of $\eta$ for which $\widehat{\alpha}$ is not defined.

By Lemma~\ref{l.dyadicCSest}, $\{n^{-(d/2+\varepsilon/4)} \alpha_n\}_{n \in \N} \in \ell^1(\N)$ for all $\varepsilon > 0$. Hence the proposition will follow from Lemma~\ref{klstool} once we prove that for all
$\varepsilon > 0$, the set of $\eta$ for which $n^{(d/2 + \varepsilon/4)} \widehat{\alpha}(\eta,n) = n^{(2d + \varepsilon)/4} \widehat{\alpha}(\eta,n)$
is unbounded is of Hausdorff dimension no more than $2d + \varepsilon$.

Recall that $d = 1-\gamma$. We consider the case where $2d < 1$, as otherwise there is nothing to prove. Choose $\varepsilon > 0$ small enough so that $2d + \varepsilon < 1$. Let $m(\eta)$ be a measurable $\Z_+$-valued function on $[0,2\pi)$. Because of
\eqref{e.vcass}, applying Lemma~\ref{SZ} with $s = 2d + \varepsilon$ yields
\begin{align*}
   \int\ \Biggl| \sum_{n=m_l(\eta)}^{2^{l+1} - 1} \alpha_n e^{i n \eta} \Biggr| \, d\nu(\eta)
&= \int\ \Biggl| \sum_{n=0}^{\widetilde{m}_l(\eta)} \alpha_{2^{l+1}-1-n} e^{-i n \eta} \Biggr| \, d \nu(\eta) \\
&\lesssim \Biggl\{ \sum_{n=2^l}^{2^{l+1}-1} (n+1)^{1 - (2d+\varepsilon)} \big|\alpha_n\big|^2
\Biggr\}^{1/2} \sqrt{\mathcal{E}_{2d+\varepsilon} (\nu)} \\
& = \Biggl\{ \sum_{n=2^l}^{2^{l+1}-1} (n+1)^{-(d+\varepsilon)} (n+1)^{\gamma} \big|\alpha_n\big|^2
\Biggr\}^{1/2} \sqrt{\mathcal{E}_{2d+\varepsilon} (\nu)} \\
&\lesssim 2^{-(d+\varepsilon) l/2}  \sqrt{\mathcal{E}_{2d+\varepsilon} (\nu)}
\end{align*}
where $m_l(\eta)=\max\{m(\eta),2^l\}$, $\widetilde{m}_l(\eta) = \min \{2^l-1, 2^{l+1} - 1 -
m(\eta)\}$, and sums with lower index greater than their upper index are to be treated as
zero. Multiplying both sides by $2^{(2d+\varepsilon) l/4}$, summing this over $l$, and
applying the triangle inequality on the left gives
$$
\int\ \Biggl| m(\eta)^{(2d+\varepsilon)/4} \sum_{n=m(\eta)}^{\infty} \alpha_n e^{i n \eta}
\Biggr| \, d \nu(\eta) \lesssim \sqrt{\mathcal{E}_{2d+\varepsilon} (\nu)}.
$$
That is, for any measurable integer-valued function $m(\eta)$, and any finite measure $\nu$ on $[0,2\pi)$,
$$
\int m(\eta)^{(2d+\varepsilon)/4} \bigl| \widehat{\alpha}(\eta, m(\eta))\bigr| \, d\nu \lesssim
\sqrt{\mathcal{E}_{2d+\varepsilon} (\nu)}.
$$
This implies that the set on which $n^{(2d+\varepsilon)/4} \widehat{\alpha}(\eta, n)$ is unbounded
must be of zero $(2d+\varepsilon)$-capacity  (i.e., it does not support a measure of finite
$(2d+\varepsilon)$-energy).

As the Hausdorff dimension of sets of zero $(2d+\varepsilon)$-capacity is less than or equal
to $2d + \varepsilon$ (see \cite[\S IV.1]{Carleson}), this completes the proof of the fact
that $S$ has Hausdorff dimension no more than $2d$.
\end{proof}

%%%%%%%%%%%%%%%%%%%%%%%%%%%%%%%%%%%%%%%%%%%%%%%%%%%%%%%%%%%%%%%%%%%%%%%%%%%%%%%%%%%%%
\subsection{Proof of Theorem~\ref{t.main} and Corollary~\ref{coro:main}}
%%%%%%%%%%%%%%%%%%%%%%%%%%%%%%%%%%%%%%%%%%%%%%%%%%%%%%%%%%%%%%%%%%%%%%%%%%%%%%%%%%%%%

We are now in a position to prove Theorem~\ref{t.main}.

\begin{proof}[Proof of Theorem~\ref{t.main}]

By Proposition~\ref{p31}, we obtain that the set
$$
S = \{ \eta \in [0,2\pi) : R_n(\eta,\beta) \text{ is unbounded for some } \beta \}
$$
has Hausdorff dimension no more than $2d = 2(1-\gamma)$.
By \eqref{fs}, $R_n \sim r_n$.
Since $r_n(\eta,0) = |\phi_n(e^{i\eta})| $ and $r_n(\eta,\pi) = |\psi_n(e^{i\eta})|$, we see that $\phi_n$ and $\psi_n$ are bounded away from the set $S$.
In view of \eqref{eq:TnzDef}, the first assertion follows.
The second assertion follows since $\mu_{{\rm sing}}$ is supported on the set $S$;
compare \cite[Corollary~10.8.4]{S05}.
\end{proof}

\begin{proof}[Proof of Corollary~\ref{coro:main}] If \eqref{eq:vcasslogdiv} holds, then, for any $\tau< \gamma$, writing $\delta = \gamma-\tau>0$, we have the following after partitioning into dyadic blocks:
\begin{align*}
\sum_{n=1}^\infty n^{\tau}|\alpha_n|^2
& =
\sum_{k=0}^\infty \sum_{n = 2^k}^{2^{k+1}-1} n^{-(\gamma-\tau)} n^{\gamma} |\alpha_n|^2 \\
& \leq
\sum_{k=0}^\infty 2^{-\delta k} \sum_{n = 2^k}^{2^{k+1}-1} n^\gamma |\alpha_n|^2 \\
& \lesssim
\sum_{k=0}^\infty 2^{-\delta k} \cdot (k+1)^B\\
& < \infty.
\end{align*}
Thus, \eqref{e.vcass} holds for all $\tau < \gamma$.
Consequently, Theorem~\ref{t.main} implies that the set of $z \in \partial \D$ for which $T_n(z)$ is unbounded has Hausdorff dimension at most $2(1-\tau)$ for all $\tau < \gamma$, so this set of $z$ has Hausdorff dimension bounded above by $2(1-\gamma)$, as desired.
\end{proof}

\end{document}